\documentclass[12pt, reqno,fleqn]{amsart}
\allowdisplaybreaks[1]
\usepackage{amsmath}
\usepackage{amssymb}
\usepackage{amsfonts}
\usepackage{verbatim}
\usepackage[usenames]{color}
\usepackage{hyperref}
\makeindex

 \newtheorem{theorem}{Theorem}[section]

\theoremstyle{definition}

\theoremstyle{remark}

\newtheorem*{claim*}{Claim}
\newtheorem{fact*}{Fact}

\newcommand{\R}{\mathbb{R}}

\newcommand{\til}{\raise.17ex\hbox{$\scriptstyle\mathtt{\sim}$}}

\newcommand\beq{\begin{equation}}

\newcommand\eeq{\end{equation}}

\newcommand\bbm{\begin{bmatrix}}
\newcommand\ebm{\end{bmatrix}}
\newcommand\bpm{\begin{pmatrix}}
\newcommand\epm{\end{pmatrix}}
\numberwithin{equation}{section}

\newlength{\Mheight}
\newlength{\cwidth}

\newcommand{\dfn}[1]{{\bf #1}\index{#1}}
\def\McCarthy{M\raise.45ex\hbox{c}Carthy}
\def\McCarthyc{M\raise.45ex\hbox{c}Carthy, }

\title[L\"owner's theorem]{The noncommutative L\"owner theorem for matrix monotone functions over operator systems}
\author{
J. E. Pascoe
}
\thanks{
$\dagger$ Partially supported by National Science Foundation Mathematical
Science Postdoctoral Research Fellowship  
DMS 1606260}
\date{\today}

\subjclass[2010]{46L52, 32A70, 30H10}

\begin{document}

\begin{abstract}
Given a function $f: (a,b) \rightarrow \mathbb{R},$
L\"owner's theorem states $f$ is monotone when extended to self-adjoint matrices via the functional calculus, if and only if $f$ extends to a self-map of the complex upper half plane. In recent years, several generalizations of L\"owner's theorem have been proven in several variables.
We use the relaxed
Agler, \McCarthyc and Young theorem on locally matrix monotone functions in several commuting variables to generalize results in the noncommutative case. Specifically, we show that a real free function defined over an operator system must analytically continue to a noncommutative upper half plane as map into another noncommutative upper half plane.
\end{abstract}
\maketitle


\section{Introduction}
Let $f: (a,b) \rightarrow \mathbb{R}.$
 L\"owner answered the question of when such a function 
is monotone when $f$ is extended to self-adjoint matrices (or even operators in general) via the functional calculus, which has found various applications.
Specifically, we say $f$ is \dfn{matrix monotone} on $(a,b) \subset \R$ if 
$$A \leq B \Rightarrow f(A) \leq f(B)$$
whenever $A$ and $B$ are self-adjoint matrices of the same size with spectrum in $(a, b),$ $f$ is being applied in the sense of the functional calculus on self-adjoint operators, and $\leq$ is interpreted to mean that the difference is positive semidefinite.

Let $\Pi$ be the upper half plane in $\mathbb{C}.$ L\"owner's theorem states the following:
\begin{theorem}[L\"owner \cite{lo34}]\label{LownerOne}
Let $f: (a,b) \to \R$ be a bounded Borel function. The function $f$ is matrix monotone if and only if $f$ is real analytic and analytically continues to the upper half plane
as a function from $\Pi \cup (a,b)$ into $\overline{\Pi}$.
\end{theorem}
For a modern treatment of L\"owner's theorem, see e.g. \cite{don74,bha97,bha07}. For applications, see e. g. \cite{Boche2004e, Jorswieck2007f, wigner, wigner2} .

Now, one could ask what kind of functions preserve inequalities of, say,
two block matrices of size $2$ by $2.$
That is, given a function $f(X_{11},X_{12},X_{21},X_{22})$ and an inequality
$$\bpm
A_{11} & A_{12} \\
A_{21} & A_{22}
\epm \leq \bpm
B_{11} & B_{12} \\
B_{21} & B_{22}
\epm,$$
when can we say
$$f(A_{11},A_{12},A_{21},A_{22}) \leq 
f(B_{11},B_{12},B_{21},B_{22})?$$
For example, is the function given by the formula
$X_{11} - X_{12}X_{22}^{-1}X_{21},$ the Schur complement, monotone in the above sense on positive block $2$ by $2$ matrices? It turns out the Schur complement is indeed monotone, which has certainly been known for some time\cite{liu99},  and can be shown via elementary arguments-- for example, its inverse appears in the formula for block inversion of a matrix as a diagonal entry. However, we are interested in an effective systematic way of classifying such functions as in L\"owner's theorem, and that is what we will establish in the noncommutative context.

\section{The noncommutative context}
We now describe the noncommutative context in which we desire to prove a generalization of L\"owner's theorem. First, we must give the appropriate generalization of the functional calculus (see \cite{vvw12} for a more thorough introduction). We also note various noncommutative generalizations to the free functional calculus of L\"owner's theorem were considered by the current author and Tully-Doyle \cite{pastd14}, and by Palfia \cite{palfia} previously, and to other functional calculi by Hansen \cite{han03} and Agler, \McCarthy, and Young \cite{amyloew}. Moreover, this work fits into a greater effort to systematize the theory of matrix inequalities \cite{heltonPositive, heltmc12, helmconvex04, helmc04, HKN}.

Let $R$ be a real topological vector space.
We define the \dfn{matrix universe over $R$,} denoted $\mathcal{M}(R)$, to be 
$$\mathcal{M}(R) = \bigcup_{n\in\mathbb{N}} M_n(\mathbb{C})\otimes_{\mathbb{R}} R,$$
where $M_n(\mathbb{C})$ denotes  the space of $n$ by $n$ matrices. We endow $\mathcal{M}(R)$ with the disjoint union topology.
Given $\mathcal{U} \subset \mathcal{M}(R),$
we use $\mathcal{U}_n$ to denote 
$\mathcal{U} \cap M_n(\mathbb{C})\otimes R$.
We define the \dfn{Hermitian matrix universe over $R$,} denoted $\mathcal{S}(R)$, to be 
$$\mathcal{S}(R) = \bigcup_{n\in\mathbb{N}} S_n(\mathbb{C})\otimes_{\mathbb{R}} R,$$
where $S_n(\mathbb{C})$ denotes  the space of $n$ by $n$ Hermitian matrices.

For a concrete example, if we take $R = S_2(\mathbb{C}),$ the $2$ by $2$ Hermitian matrices over $\mathbb{C},$ $\mathcal{M}(S_2(\mathbb{C}))$ consists of all block $2$ by $2$ matrices and $\mathcal{S}(S_2(\mathbb{C}))$ consists of all block $2$ by $2$ Hermitian matrices. For another example, taking $R = \mathbb{R}^2,$ the set $\mathcal{M}(\mathbb{R}^2)$ consists of all pairs of same-sized matrices and $\mathcal{S}(\mathbb{R}^2)$ consists of all pairs of same-sized Hermitian matrices.

We define a {\bf domain} $D \subset \mathcal{M}(R)$
to satisfy the following two axioms:
\begin{enumerate}
	\item $X \oplus Y \in D \Leftrightarrow X, Y \in D$
	\item $X \in D_n \Rightarrow U^*XU \in D$ for all $n$ by $n$ unitaries 
	$U$ over $\mathbb{C}.$
\end{enumerate}
Let $D \subset \mathcal{M}(R_1)$ be a free domain.
We say a function $f:D \rightarrow \mathcal{M}(R_2)$ is a {\bf free function} if
\begin{enumerate}
	\item $f|_{D_n}$ maps into $\mathcal{M}(R_2)_n$
	\item $f(X\oplus Y) = f(X) \oplus f(Y),$
	\item $S^{-1}f(X)S = f(S^{-1}XS)$ for all $n$ by $n$ invertible matrices $S$
	 over $\mathbb{C}$ such that $X, S^{-1} XS \in D_n.$
\end{enumerate}
We note any noncommutative rational expression gives a free function on its domain of definition. For example, the Schur complement, $X_{11} - X_{12}X_{22}^{-1}X_{21},$ gives a free function on the subset $D \subseteq \mathcal{M}(S_2(\mathbb{C}))$ where $X_{22}^{-1}$ is defined. For another example, the matrix geometric mean,
$X_1^{1/2}(X_1^{-1/2}X_2X_1^{-1/2})^{1/2} X_1^{1/2}$
defines a free function on the subset $D \subseteq \mathcal{S}(\mathbb{R}^2)$ of all pairs of positive definite matrices.

If $R$ is a real operator system, that is, $R$  is a real subspace containing $1$ in a $C^*$ algebra of self-adjoint elements, for each $n$ there is a natural ordering on $S_n(\mathbb{C})\otimes R,$ since matrices over $R$ are themselves elements of a larger $C^*$-algebra. That is,  given $A, B \in S_n(\mathbb{C})\otimes R,$ we say $A \leq B$ if $B - A$ is positive semidefinite as an element of  $S_n(\mathbb{C})\otimes R.$

Accordingly, given $R_1$ and $R_2$ real operator systems and a domain $D \subseteq \mathcal{S}(R_1)$, we say a free function $f: D \rightarrow \mathcal{S}(R_2)$ is \dfn{matrix monotone} if 
$A \leq B \Rightarrow f(A) \leq f(B)$
whenever $A$ and $B$ have the same size.

Define 
$\Pi(R) = \{\text{Im } X > 0\}$
where $A > B$ if the difference is strictly positive definite, that is, it is self-adjoint and its spectrum is a subset of $(0,\infty)$, and $\text{Im } X = (X - X^*)/2i$.

We show the following theorem.
\begin{theorem}[Noncommutative L\"owner theorem over operator systems]
Let $R_1$ and $R_2$ be closed real operator systems.
Let $D \subseteq \mathcal{S}(R_1)$ be a free domain.
Suppose each $D_n$ is convex and open as a subset of $S_n(\mathbb{C})\otimes R$.
A function $f: D \rightarrow \mathcal{S}(R_2)$ is matrix monotone
if and only if
$f$ extends to a continuous free function $F: \Pi(R_1) \cup D \rightarrow \overline{\Pi(R_2)}.$
\end{theorem}
We note that such a function must be analytic on each $\Pi(R_1)_n$ due to the draconian nature of free functions. See \cite{vvw12}. We also point out that the case where $R_1 = \mathbb{R}^d$ as a diagonal algebra and $R_2 = \mathbb{R}$ was explored in \cite{pastd14, palfia}, and that the current work simplifies the proof of the main result of those works if we are willing to use the commutative L\"owner theorem from \cite{amyloew} as a black box.  Moreover, if we are given a rational expression, such as the Schur complement, on a nice finite dimensional operator system, such as a matrix algebra,  one can apply the algorithms in \cite{HKN} which make the rational convex Positivstellensatz
\cite{pascoePosSS} effective to check that a function is matrix monotone in our sense. 

Finally, we should comment that the setting of operator systems is equivalent to defining an Archimedian matrix ordering
on  $\mathcal{S}(R)$,  where $R$ is an abstract real vector space, by the Choi-Effros Theorem \cite{effrosChoi}.
That is, we might have alternatively defined an ordering on $\mathcal{S}(R)$ using any
proper closed Archimedian matrix convex cone, but the result is the same.

Before we arrive at the proof of our Theorem, we should revisit our Schur complement. Our domain $D \subset \mathcal{S}(S_2(\mathbb{C}))$ is the set of positive definite block $2$ by $2$ matrices  upon which  our function, defined by the formula
$$f\bpm
X_{11} & X_{12} \\
X_{21} & X_{22}
\epm = X_{11} - X_{12}X_{22}^{-1}X_{21},$$
is a free function $f: D \rightarrow \mathcal{S}(\mathbb{R}).$
According to our Theorem, $f$ will be matrix monotone if and only if $f$ extends to a continuous free function from
$D \cup  \Pi(S_2(\mathbb{C}))$ to $\overline{\Pi(\mathbb{R})}.$
It is clear that extension of $f$ to the new domain must still be given by the same formula as before.
Either using the algorithms in \cite{HKN, pascoePosSS} or by brute force, one can see that
$$\textrm{Im } f =
\bpm 1 \\ (X_{22}^*)^{-1} X_{12}^* \epm^*
\left[ \textrm{Im } \bpm
X_{11} & X_{12} \\
X_{21} & X_{22}
\epm\right] \bpm 1 \\ (X_{22}^*)^{-1}X_{12}^* \epm$$
which is manifestly positive definite whenever  $\textrm{Im } \bpm
X_{11} & X_{12} \\
X_{21} & X_{22}
\epm$ is positive definite-- that is $f$ maps $\Pi(S_2(\mathbb{C}))$ to $\Pi(\mathbb{R}).$ That is, our Theorem now implies that the Schur complement is matrix monotone.

Another example of a matrix monotone function, is the matrix geometric mean and various generalizations, see \cite{Lawson2011, Bhatia2012}. In the two parameter case it is not immediately clear to the author how to show the 
function
$X_1^{1/2}(X_1^{-1/2}X_2X_1^{-1/2})^{1/2} X_1^{1/2}$
continues to a map from $\Pi(\mathbb{R}^2)$ to
$\overline{\Pi(\mathbb{R})}$ without going through the generalization of L\"owner's theorem.

\section{The proof of the main result}
$(\Rightarrow)$
The proof will go by viewing, for each $n,$ $f|_{D_n}$ as a matrix monotone function in several commuting variables in the sense of Agler, \McCarthy, and Young.

Agler, \McCarthy, and Young extended L\"owner's theorem to several commuting
variables for the class of locally matrix monotone functions \cite{amyloew}.
Subsequently, it was generalized to remove some technical assumptions by the author in \cite{pascoemollifier}.
Let $E$ be an open subset of $\mathbb{R}^d.$
Let $CSAM^d_n(E)$ denote the $d$-tuples of commuting self-adjoint matrices of size $n$ with joint spectrum contained in $E$.
 We say that a function $f: E \rightarrow \mathbb{R}$ is \dfn{locally matrix monotone} if for any $C^1$
path $\gamma: (-1,1) \rightarrow CSAM^d_n(E)$ such that $\gamma'(0)_i > 0$, there exists
an $\epsilon >0$ such that for all $-\epsilon < t_1 < t_2 < \epsilon,$
$g(\gamma(t_1)) \leq g(\gamma(t_2)).$

 We recall the following theorem.
\begin{theorem}[Agler, \McCarthy, and Young \cite{amyloew}, Pascoe \cite{pascoemollifier}] \label{amyold}
	Let $E$ be an open subset of $\mathbb{R}^d.$
	Let $g: E \rightarrow \mathbb{R}$ be a locally matrix monotone function. Then $g$ is analytic, and $g$ extends to a (unique) continuous function on $\Pi^d \cup E$ which maps into $\overline{\Pi}$ which is analytic on $\Pi^d.$
\end{theorem}
We note that the original formulation of Agler, \McCarthy, and Young applied only to $C^1$ functions $g$, and via an argument using mollifiers it was shown that the theorem holds for arbitrary functions.

We note that is sufficient to show that on each nonempty $D_n$, our function $f$ analytically continues to $\Pi(R_1)_n$ taking values in $\overline{\Pi(R_2)_n}$. It is an elementary, but perhaps somewhat involved, exercise to show that the induced extension of $f$ will be a free function on $\Pi(R_1).$ Namely, the edge-of-the-wedge theorem will ensure that the extension of $f$ actually analytically continues through each $D_n$ as a function on $\Pi(R_1)_n \cup D_n$ and the rest of the properties will follow by analytic continuation.

Now, we note that it is sufficient to show that for every (completely) positive unital linear functional $l: M(R_2)_n \rightarrow \mathbb{C}$ that $f_{n,l} = l \circ f|_{D_n}$
extends analytically to $D_n \cup \Pi(R_1)_n$ as a map taking values with positive imaginary part.
This is obvious when $R_2$ is finite dimensional, and an exercise in functional analysis otherwise.

Fix $P \in D_n.$
Let $K_1, \ldots K_m > 0$ be positive elements of $\mathcal{S}(R_1)_n.$
Let $C$ be the cone generated by $K_1, \ldots, K_m$
and let $S$ be the span of $K_1, \ldots, K_m$.
We will show that $f_{l,n}$ uniquely analytically continues to
$P + S + iC.$  Taking larger and larger sets of $K_i$ will give an analytic continuation $f_{n,l}$ to the whole of $\Pi(R_1)_n.$ That is, the sets $P +S+ iC$ exhaust $\Pi(R_1)_n.$

Define the function $g(h) = f_l(P + \sum_i h_i K_i))$.
Now $g(h)$ is a locally matrix monotone function in the sense of Theorem \ref{amyold} as a function on $\mathbb{R}^m$ which induces the unique analytic continuation of $f_{n,l}$ to 
the desired space taking values in $\overline{\Pi}.$ So, we are done.

$(\Leftarrow)$ The converse direction is easy and follows from a computation of the derivative for directions pointing into the upper half plane. See \cite[Lemma 4.8]{pastd14} where the details are essentially the same.

\bibliography{references}
\bibliographystyle{plain}

\end{document}